\documentclass[11pt]{book}

 \usepackage{a4}
\usepackage{graphicx}
\usepackage{amsfonts}
\usepackage{natbib}
\usepackage{algorithm}
\usepackage{algorithmic}
\usepackage{parskip}

\usepackage{latexsym}

\newtheorem{theorem}{Theorem}[section]
\newtheorem{lemma}{Lemma}[section]
\newtheorem{corollary}{Corollary}[section]
\newtheorem{remark}{Remark}[section]
\newtheorem{example}{Example}[section]
\newtheorem{definition}{Definition}[section]

 \newcommand{\vvert} {| \hskip-0.15em |  \hskip-0.15em |}
\newcommand{\R}{\mathbb{R}}
\newcommand{\Nat}{\mathbb{N}}

\newcommand{\PP} {{  \rm I\hskip-0.22em P}}

\newcommand{\EE} {{\rm I\hskip-0.48em E}}

\usepackage{graphicx}
\usepackage{parskip}

\usepackage{amsfonts}
\usepackage{natbib}
\usepackage{algorithm}
\usepackage{algorithmic}

\usepackage{latexsym}

\newtheorem{proposition}{Proposition}[section]

\usepackage{natbib}

\title{Statistical Theory for High-Dimensional Models
\vskip .5in
\centerline{Lecture Notes}}

\date{September 2014}
\author{Sara van de Geer}

\begin{document}

\maketitle
\tableofcontents
\newpage

\chapter{The Lasso and variations}

{\bf Abstract} {\it We present oracle inequalities for the prediction error of the Lasso and
square-root Lasso and briefly describe the scaled Lasso.}

\section{The model}

Let $Y \in \R^n$ be an $n$-vector of real-valued observations and let $X$
be a given $n \times p $ design matrix.
We let
$$\EE Y := f_0 . $$
We assume $X$ to be fixed, i.e., we consider the case of fixed design. The entries of the vector
$f_0$ are thus the (conditional) expectation of $Y$ given $X$. 
We let $ \epsilon := Y-f^0 $
be the noise term. 

We assume that $X$ has rank $n$ so that there is a solution $\beta^0$ of the equation $f^0 = X \beta^0 $.
We may then take e.g. the basis pursuit solution (\cite{chen1998atomic})
$$ \beta^0 := \arg \min \{ \| \beta \|_1 : \ X \beta = f^0 \}  . $$

\section{Notation}

For a vector $v \in \R^n$ we use the notation $\| v \|_n:= \sqrt { v^T v /n} $. 
Write the (normalized) Gram matrix as $\hat \Sigma := X^T X / n $. Thus $\| X \beta \|_n^2
= \beta^T \hat \Sigma \beta$, $\beta \in \R^p$.

For a vector $\beta$ we denote its $\ell_1$-norm as
$\| \beta \|_1$.
The dual norm of $\| \cdot \|_1$ is the $\ell_{\infty} $-norm $\| \cdot \|_{\infty}$.
The {\it dual norm inequality} says that for any two vectors $w$ and $\beta$ 
$$| w^T \beta | \le  \| w \|_{\infty} \| \beta \|_1. $$

Let $S \subset \{ 1 , \ldots , p \}$ be an index set. 
We use the notation
$$\beta_{j,S} := \beta_j {\rm l} \{ j \in S \} , \ j=1 , \ldots , p. $$
Thus $\beta_S$ is a $p$-vector with entries equal to zero at the indexes $j \notin S$.
We will sometimes identify $\beta_S$ with the vector
$\{ \beta_j \}_{j \in S} \in \R^{|S|}$.  We let $S_{\beta} := \{ j : \ \beta_j \not= 0 \} $ be the active
set of the vector $\beta$. 

\begin{definition} (\cite{vandeGeer:07a}, \cite{BvdG2011})
For constant $L>0$ and an index set $S$ the compatibility constant  is
$$ \hat \phi^2 (L,S) := \min \biggl \{ |S|  \| X \beta_S - X \beta_{S^c} \|_n^2 : \ \| \beta_S \|_1 = 1 , \ \| \beta_{S^c } \|_1 \le L \biggr \}. $$
\end{definition} 

\section{The Lasso}

The Lasso estimator (\cite{tibs96}) $\hat \beta$ is defined as
$$ \hat \beta := \arg \min_{\beta \in \R^p} \biggl \{ \| Y - X \beta \|_n^2 +
2 \lambda \| \beta \|_1 \biggr \} . $$

This estimator satisfies the {\it Karush-Kuhn-Tucker conditions} or KKT-conditions
which say that
\begin{equation} \label{KKT.equation}
X^T (Y- X \hat \beta ) /n = \lambda \hat z 
\end{equation}
where $\hat z $ is a $p$-dimensional vector with $\| \hat z \|_{\infty} \le 1 $ and
with $\hat z_j = {\rm sign} (\hat \beta_j)$ if $\hat \beta_j \not= 0 $. The KKT-conditions
follow from sub-differential calculus which defines the sub-differential of
the absolute value function $x \mapsto |x|$ as 
$$\{ {\rm sign} (x) \} \{ x \not= 0 \} + [-1,1] \{ x=0 \} . $$
As a consequence we have the {\it KKT}-inequality: for any $\beta \in \R^p$
$$ (\beta - \hat \beta )^T X^T (Y - X \hat \beta) /n \le \lambda \| \beta \|_1 -
  \lambda \| \hat \beta\|_1 .$$
  As we will see in our proofs this inequality is useful in conjunction with the {\it three beta's} layout
  $$
   2 (\hat \beta - \beta )^T \hat \Sigma  (\hat \beta - \beta^0) = \| X ( \hat \beta - \beta^0 ) \|_n^2  - \| X (\beta - \beta^0 ) \|_n^2 +
   \| X (\hat  \beta -  \beta ) \|_n^2  .$$
  
  Another important inequality will be the {\it convex conjugate} inequality: for any $a,b \in \R$
  $$ 2ab \le a^2 + b^2 . $$ 
  
  We will also often use the $\ell_1$-triangle trick: suppose for some $\beta, \tilde \beta$ and constant $c$,
  $$\| \tilde \beta \|_1 \le \|  \beta \|_1 + c , $$
  then
  $$ \|\tilde \beta_{S^c} \|_1  \le \| \tilde \beta_S - \beta \|_1 +c , $$
  where $S= S_{\beta} $ is the active set of $\beta$.

The random vector $\epsilon^T X$
occurring below has mean zero. To control its $\ell_{\infty}$-norm we will use empirical process theory (see 
Section \ref{pmaximum.section}). 

\begin{theorem}\label{Lasso.theorem} (\cite{koltchinskii2011nuclear})
Let $\lambda_{\epsilon}$ satsify
$$\lambda_{\epsilon} \ge  \| \epsilon^T X  \|_{\infty} /n .  $$ 
Define  for $\lambda >   \lambda_{\epsilon}  $
$$L:= { \lambda + \lambda_{\epsilon} \over \lambda- \lambda_{\epsilon}} . $$
Then
$$ \| X ( \hat \beta - \beta^0 ) \|_n^2  \le \min_{S} \biggl \{ \min_{\beta \in \R^p,\  S_{\beta}   =S  }  \| X( \beta - \beta^0) \|_n^2 +  (\lambda + \lambda_{\epsilon} )^2 |S | / \hat \phi^2 (L,S) \biggr \} .  $$
\end{theorem}

 {\bf Proof.} 
   Fix some $\beta \in \R^p$  and let $S := \{ j: \ \beta_j \not= 0 \} $ be its active set.
  If 
  $$( \hat \beta - \beta)^T \hat \Sigma ( \hat \beta -\beta^0 ) \le 0$$
  we find from the three beta's layout
  $$\| X (\hat \beta - \beta^0)\|_n^2  $$ $$=
\| X (\beta - \beta^0) \|_n^2
     - \| X (\beta - \hat \beta)\|_n^2  + 2 ( \hat \beta - \beta)^T \hat \Sigma  ( \hat \beta - \beta^0 ) $$ $$\le  
  \| X (\beta - \beta^0)\|_n^2  . $$
  Hence then we are done.
  
  Suppose now that
  $$( \hat \beta - \beta^0)^T \hat \Sigma ( \hat \beta -\beta )  \ge 0  .$$
  By the KKT-inequality  
  $$(\beta - \hat \beta )^T X^T (Y - X \hat \beta) /n \le \lambda \| \beta \|_1 -
  \lambda \| \hat \beta\|_1 .$$
    As $Y= X \beta^0 + \epsilon$
  $$(\hat \beta -\beta )^T \hat \Sigma ( \hat \beta - \beta^0  ) + \lambda \| \hat \beta \|_1  \le \epsilon^T X(\hat \beta - \beta)/n + \lambda \| \beta \|_1  . $$
  By the dual norm inequality
  $$| \epsilon^T X(\hat \beta - \beta)|/n \le (\| \epsilon^T X \|_{\infty}/n)  \| \hat \beta - \beta \|_1 \le
  \lambda_{\epsilon}  \| \hat \beta - \beta \|_1
  . $$
  Thus
  $$ (\hat \beta -\beta )^T \hat \Sigma ( \hat \beta - \beta^0  )  + \lambda \| \hat \beta \|_1 \le 
  \lambda_{\epsilon} \| \hat \beta - \beta \|_1  + \lambda \| \beta\|_1 . $$
  By the $\ell_1$-triangle trick this implies
  \begin{equation}\label{tea.equation}
  (\hat \beta -\beta )^T \hat \Sigma ( \hat \beta - \beta^0  )  + (\lambda - \lambda_{\epsilon})  \| \hat \beta_{S^c}\|_1 \le 
   (\lambda +  \lambda_{\epsilon} ) \| \hat \beta_S - \beta \|_1   . 
   \end{equation} 
 Since $  ( \hat \beta - \beta)^T \hat \Sigma ( \hat \beta -\beta^0 )  \ge 0$ this gives
 $$\| \hat \beta_{S^c} \|_1  \le L \|\hat \beta_S - \beta\|_1 . $$
 By the definition of the compatibility constant $\hat \phi^2 (L,S)$ we then have
 \begin{equation}\label{coffee.equation}
\| \hat \beta_S - \beta \|_1  \le \sqrt {|S|}  \| X ( \hat \beta - \beta ) \|_n / \hat \phi(L,S) .
 \end{equation}
 Continue with inequality (\ref{tea.equation}) and apply the convex conjugate inequality
 $$ (\hat \beta -\beta )^T \hat \Sigma ( \hat \beta - \beta^0  ) 
   + (\lambda - \lambda_{\epsilon} )  \| \hat \beta_{S^c}\|_1 $$ $$ \le 
   (\lambda + \lambda_{\epsilon} )\sqrt {|S|}   \| X ( \hat \beta - \beta ) \|_n /\hat \phi (L, S)  $$
   $$ \le {1 \over 2} { |S|  (\lambda + \lambda_{\epsilon}  )^2 \over  \hat \phi^2(L,S) }   + {1 \over 2} \| X (\hat \beta - \beta ) \|_n^2
     . $$
   Since by the three beta's layout
   $$2 ( \hat \beta - \beta^0)^T \hat \Sigma ( \hat \beta -\beta )  
   = \| X ( \hat \beta - \beta^0 ) \|_n^2  - \| X (\beta - \beta^0 ) \|_n^2 +
   \| X ( \hat \beta - \beta ) \|_n^2  ,$$
   we obtain
   $$\| X (\hat \beta - \beta^0 ) \|_n^2 + 2(\lambda - \lambda_{\epsilon} ) \| \hat \beta_{S^c}\|_1   
\le \| X (\beta - \beta^0 ) \|_n^2 +  |S|(\lambda + \lambda_{\epsilon}  )^2 / \hat \phi^2 (L,S)  . $$
  \hfill $\sqcup \mkern -12mu \sqcap$
  
  We will now show that if one increases the constant $L$ in the compatibility constant, one
  can prove a bound for the $\ell_1$-estimation error. 
  
   \begin{theorem}\label{ell1.theorem} 
Let $\lambda_{\epsilon}$ satsify
$$\lambda_{\epsilon} \ge  \| \epsilon^T X  \|_{\infty} /n .  $$ 
Let $0 \le \delta < 1$ be arbitrary and define  for $\lambda >   \lambda_{\epsilon}  $
$$L:=  { \lambda + \lambda_{\epsilon} + \delta ( \lambda - \lambda_{\epsilon} )\over (1- \delta) (\lambda- \lambda_{\epsilon}) } . $$
Then 
$$ 2\delta  (\lambda - \lambda_{\epsilon} ) \| \hat \beta - \beta^0 \|_1+  \| X ( \hat \beta - \beta^0 ) \|_n^2  $$ $$
\le  \min_{S} \biggl \{ \min_{\beta \in \R^p, \ S_{\beta} = S } \biggl [2\delta (\lambda - \lambda_{\epsilon} ) \| \beta - \beta^0 \|_1 + 
 \| X( \beta - \beta^0) \|_n^2  \biggr ]  $$ $$
+  \biggl [ \lambda + \lambda_{\epsilon} +  \delta (\lambda -\lambda_{\epsilon}) \biggr ]^2 { |S | \over  \hat \phi^2 (L,S) }\biggr \} .  $$
\end{theorem}

{\bf Proof.} We follow the same line of reasoning as in the proof of Theorem \ref{Lasso.theorem}.
Let $\beta \in \R^p$ and $S:= \{ j :\ \beta_j \not= 0 \}$. 
    If 
  $$\| X (\hat \beta - \beta) \|_n^2  \le -\delta (\lambda - \lambda_{\epsilon} ) \|  \hat \beta - \beta \|_1$$  we find from the three beta's layout
  $$2\delta (\lambda - \lambda_{\epsilon})\|  \hat \beta - \beta \|_1  +\| X ( \hat \beta - \beta^0 )\|_n^2  
  $$ $$ =2\delta (\lambda - \lambda_{\epsilon} ) \| \hat \beta - \beta \|_1  + \| X (\beta - \beta^0 ) \|_n^2 
      - \| X ( \beta - \hat \beta ) \|_n^2  + 2 ( \hat \beta - \beta)^T \hat \Sigma  ( \hat \beta - \beta^0 ) $$ $$
    \le  \| X (\beta - \beta^0 ) \|_n^2 $$
   and we are done.
  
  Suppose now that
  $$(\hat \beta - \beta)^T \hat \Sigma ( \hat \beta -\beta^0 ) \ge -\delta (\lambda - \lambda_{\epsilon} ) \|  \hat \beta - \beta \|_1
  .$$
  By the KKT-inequality we have
  $$(\hat \beta - \beta )^T \hat \Sigma (\hat \beta -\beta^0 )  + \lambda \| \hat \beta \|_1  \le \epsilon^T X(\hat \beta - \beta)/n + \lambda \|\beta \|_1 . $$
  By the dual norm inequality
  $$| \epsilon^T X(\hat \beta - \beta)| /n   \le \lambda_{\epsilon} \| \hat \beta - \beta \|_1 
    . $$
  Thus
  $$ (\hat \beta - \beta )^T\hat \Sigma (\hat \beta -\beta^0 )  + \lambda \| \hat \beta\|_1 \le 
   \lambda_{\epsilon} \| \hat \beta - \beta \|_1 + \lambda \| \beta\|_1 . $$
  By  the $\ell_1$-triangle trick this implies
  \begin{equation}\label{tea1.equation}
  (\hat \beta - \beta )^T\hat \Sigma (\hat \beta -\beta^0 )  + (\lambda - \lambda_{\epsilon} )  \| \hat \beta_{S^c} \|_1 \le 
   (\lambda +  \lambda_{\epsilon} ) \| \hat \beta_S - \beta \|_1  . 
   \end{equation} 
 Since $  (\hat \beta - \beta )^T \hat \Sigma (\hat \beta -\beta^0 ) \ge - \delta (\lambda - \lambda_{\epsilon} ) \| \hat \beta - \beta \|_1$ this gives
 $$(1-\delta) ( \lambda - \lambda_{\epsilon}) \|  \hat \beta_{S^c} \|_1 \le (\lambda + \lambda_{\epsilon}+
 \delta ( \lambda - \lambda_{\epsilon} ) ) \| \hat \beta_S - \beta \|_1 $$
 or
 $$ \|  \hat \beta_{S^c} \|_1 \le  L \| \hat \beta_S - \beta \|_1.  $$
 But then
 \begin{equation}\label{coffee1.equation}
 \| \hat \beta_S - \beta \|_1  \le \sqrt {|S|}   \| X ( \hat \beta - \beta ) \|_n / \hat \phi(L,S) .
 \end{equation}
 Continue with inequality (\ref{tea1.equation}) and apply the convex conjugate inequality:
 $$(\hat \beta - \beta )\hat \Sigma  (\hat \beta -\beta^0 )  + (\lambda - \lambda_{\epsilon})  \| \hat \beta_{S^c}\|_1 + \delta ( \lambda - \lambda_{\epsilon}) 
 \| \hat \beta_S - \beta \|_1 $$ $$ \le 
   [  \lambda + \lambda_{\epsilon} +  \delta (\lambda -\lambda_{\epsilon})  ] \sqrt {|S|}  \| X ( \hat \beta - \beta ) \|_n  /
   \hat \phi (L,S) $$
   $$ \le {1 \over 2} \biggl [  \lambda + \lambda_{\epsilon}  + \delta (\lambda -\lambda_{\epsilon}) \biggr ]^2 { |S|  \over  \hat \phi^2 (L,S) }  + {1 \over 2} \| X ( \hat \beta - \beta ) \|_n^2 . $$
   Since by the three beta's layout
   $$
   2 (\hat \beta - \beta )^T \hat \Sigma ( \hat \beta - \beta^0 )= \| X ( \hat \beta - \beta^0 ) \|_n^2  - \| X (\beta - \beta^0 ) \|_n^2 +
   \| X ( \hat \beta - \beta ) \|_n^2  ,$$
   we obtain
   $$\| X (\hat \beta - \beta^0 ) \|_n^2 + 2(\lambda - \lambda_{\epsilon} )  \| \hat \beta_{S^c}\|_1  +  2 \delta ( \lambda - \lambda_{\epsilon})   \|  \hat \beta_S - \beta\|_1 $$ $$
\le \| X (\beta - \beta^0 ) \|_n^2 + \biggr [ \lambda + \lambda_{\epsilon} + \delta (\lambda -\lambda_{\epsilon})  \biggr ]^2 |S|  /\hat \phi^2 (L,S)  . $$
  \hfill $\sqcup \mkern -12mu \sqcap$
  
  The result of Theorem \ref{ell1.theorem} leads to a trade-off between the approximation error
  $\| X ( \beta - \beta^0) \|_n^2$, the $\ell_1$-error $\| \beta - \beta^0 \|_1$ and the sparseness\footnote{or non-sparseness actually} 
  $S_{\beta}$ (or rather the {\it effective sparseness}  $S_{\beta}/ \hat \phi^2 (L,S_{\beta} ) $).
  To study this let us consider the oracle $\beta^*$ which trades off approximation error and
  (effective) sparsity but is meanwhile restricted to have an $\ell_1$-norm at least as large as that of $\beta^0$.
  
  \begin{lemma} \label{ell1-restricted.lemma} Let for some $\lambda_*>0$ the vector $\beta^*$ be defined as
  $$ \beta^* := \arg \min \biggl \{ \| X ( \beta - \beta^0 )\|_n^2 + \lambda_*^2 |S_{\beta} | / \hat \phi^2 (L,S_{\beta}) :\
  \| \beta \|_1 \ge \| \beta^0 \|_1 \biggr \} . $$
  Let $S_* := \{ j:\ \beta_j^* \not= 0 \}$ be the active set of $\beta^*$. 
  Then
  $$ \lambda_* \| \beta^* - \beta^0 \|_1 \le \| X ( \beta^* - \beta^0 ) \|_n^2 + {\lambda_*^2 | S_* | \over \hat \phi^2 (1,S_*) } 
  . $$
  \end{lemma}
  
  {\bf Proof.} Let $S_* := \{ j: \ \beta_j^* \not= 0 \} $. Since $\| \beta^0 \|_1 \le \| \beta^* \|_1 $ we know
  by the $\ell_1$-triangle trick
  $$\| \beta_{S_*^c}^0 \|_1  \le \| \beta^* - \beta_{S_*}^0 \|_1 . $$
  Hence by the definition of the compatibility constant and by the convex conjugate inequality
  $$ \lambda_* \| \beta^* - \beta^0 \|_1 \le 2 \lambda_* \| \beta^* - \beta_{S_*}^0 \|_1 \le 
  { 2 \lambda_* \| X ( \beta^* - \beta^0 ) \|_n \over  \hat \phi (1, S_*) } \le
  \| X ( \beta^* - \beta^0 ) \|_n^2 + {\lambda_*^2 | S_* | \over \hat \phi^2 (1,S) } . $$
  \hfill $\sqcup \mkern -12mu \sqcap$
  
  From Lemma \ref{ell1-restricted.lemma} we see that an $\ell_1$-restricted oracle $\beta^*$ that trades off approximation error
  and sparseness is also going to be close in $\ell_1$-norm. We have the following corollary for the bound of
  Theorem \ref{ell1.theorem}.
  
  \begin{corollary}\label{ell1-restricted.corollary} 
  Let 
$$\lambda_{\epsilon} \ge  \| \epsilon^T X  \|_{\infty} /n .  $$ 
Let $0 \le \delta < 1$ be arbitrary and define  for $\lambda >   \lambda_{\epsilon}  $
$$L:=  { \lambda + \lambda_{\epsilon} + \delta ( \lambda - \lambda_{\epsilon} )\over (1- \delta) (\lambda- \lambda_{\epsilon}) } . $$
Let the vector $\beta^*$ with active set $S_*$ be defined as in Lemma 
  \ref{ell1-restricted.lemma} with $\lambda_* :=  \lambda+ \lambda_{\epsilon} +
  \delta ( \lambda - \lambda_{\epsilon} )$.
We have 
$$2\delta  (\lambda - \lambda_{\epsilon} ) \| \hat \beta - \beta^0 \|_1 + \| X ( \hat \beta - \beta^0 ) \|_n^2   $$ $$
    \le  {\lambda + \lambda_{\epsilon} + 3\delta(\lambda - \lambda_{\epsilon} )  \over \lambda + \lambda_{\epsilon}+ \delta (\lambda -\lambda_{\epsilon})} 
 \| X ( \beta^* - \beta^0 )\|_n^2 $$ $$ + \biggl ( \lambda + \lambda_{\epsilon} + 3 \delta (\lambda -\lambda_{\epsilon}) \biggr )
 \biggl ( \lambda + \lambda_{\epsilon} +\delta (\lambda - \lambda_{\epsilon}) \biggr ) 
{  |S_* |  /   \hat \phi^2 (L,S_*) }  .$$  
    
  \end{corollary}

\section{The square-root Lasso}

In the previous section we required that the tuning parameter $\lambda$ is chosen at least
as large as the {\it noise level} $\lambda_{\epsilon}$ where $\lambda_{\epsilon}$ is a bound
for $\| \epsilon^T X\|_{\infty} / n$. Clearly, if for example the entries in $\epsilon$ are i.i.d.\ with variance 
$\sigma^2$, the choice of $\lambda$ will depend on the standard deviation $\sigma$ which will
usually be unknown in practice. To avoid this problem we 
consider the square-root Lasso (\cite{belloni2011square})
$$ \hat \beta := \arg \min_{\beta \in \R^p} \biggl \{ \| Y - X \beta \|_n +
 \lambda \| \beta \|_1 \biggr \} . $$
  
 The square-root Lasso  $\hat \beta$ satisfies the KKT-conditions
 \begin{equation} \label{square-root-KKT.equation}
{ X^T ( Y- X \hat \beta ) /n \over \| Y - X \hat \beta \|_n} = \lambda \hat z 
\end{equation}
where $\| \hat z \|_{\infty} \le 1$ and $\hat z_j = {\rm sign}
(\hat \beta_j)$ if $\hat \beta_j \not= 0 $. 
Defining the residuals $\hat \epsilon := Y - X \hat \beta $ we can write this as
$$X^T ( Y- X \hat \beta ) /n  = \lambda \| \hat \epsilon \|_n \hat z .$$

\begin{proposition}\label{square-rootLasso.proposition} Let $\hat \lambda_0$ satisfy
$$\hat \lambda_0 \| \hat \epsilon \|_n  \ge  \| \epsilon^T X  \|_{\infty} /n .  $$ 
Define  for $\lambda >  \hat \lambda_0  $
$$\hat L:= { \lambda + \hat \lambda_0 \over \lambda- \hat \lambda_0} . $$
Then
$$ \| X ( \hat \beta - \beta^0 ) \|_n^2  \le \min_{S} \min_{\beta \in \R^p,\  S_{\beta}   =S  } \biggl \{ \| X( \beta - \beta^0) \|_n^2 +  (\lambda + \hat \lambda_0 )^2 \| \hat \epsilon \|_n^2 { |S | \over  \hat \phi^2 (\hat L,S) }\biggr \} .  $$
\end{proposition}

{\bf Proof.} The estimator $\hat \beta$  
satisfies the KKT-conditions (\ref{square-root-KKT.equation}) which are exactly
the KKT-conditions (\ref{KKT.equation}) but with $\lambda$ replaced by
$\lambda \| \hat \epsilon \|_n$.  This means we can recycle the proof of Theorem \ref{Lasso.theorem}.
\hfill $\sqcup \mkern -12mu \sqcap$

Proposition \ref{square-rootLasso.proposition} is not very useful as such
because it depends on $\| \hat \epsilon \|_n$. We therefore provide bounds for $\| \hat \epsilon \|_n$.

\begin{lemma} \label{square-root-bounds.lemma} We have $\| \hat \epsilon \|_n \le \| \epsilon \|_n + \lambda \| \beta^0 \|_1 $. If for a constant $\lambda_0$ satisfying $\lambda_0  \| \epsilon \|_n \ge\| \epsilon^T X \|_{\infty} / n $ the tuning parameter $\lambda$ has
$\lambda(1- \eta^2) > 2 \lambda_0$ for some $\eta > 0$ and 
$$ \| \beta^0 \|_1 \le {\| \epsilon \|_n \over 4 \lambda_0 } \biggl ( 1 - \eta^2- { 2 \lambda_0 \over \lambda }   \biggr ) . $$
then $\| \hat \epsilon \|_n^2 \ge \eta^2 \| \epsilon \|_n^2 $.

\end{lemma}

{\bf Proof.} 
Since $\hat \beta $ minimizes $\| Y - X \beta \|_n + \lambda \| \beta \|_1 $ we have
$$\| \hat \epsilon \|_n = \| Y - X \hat \beta \|_n \le \| Y - X \beta^0 \|_n + \lambda \| \beta^0 \|_1 -
\lambda \| \hat \beta \|_1$$ $$ \le  \| Y - X \beta^0 \|_n + \lambda \| \beta^0 \|_1 = \| \epsilon \|_n + 
\lambda \| \beta^0 \|_1 . $$
Moreover
$\| \hat \beta \|_1 \le { (\| \epsilon \|_n - \| \hat \epsilon \|_n)  / \lambda} + \| \beta^0 \|_1 \le
\| \epsilon \|_n / \lambda + \| \beta^0 \|_1 $.
Hence
$$ \| \hat \epsilon \|_n^2 = \| \epsilon - X ( \hat \beta - \beta^0 ) \|_n^2 \ge
\| \epsilon \|_n^2 - 2 \lambda_0 \| \epsilon \|_n \| \hat \beta - \beta^0 \|_1  $$
$$ \ge \| \epsilon \|_n^2 - 2 \lambda_0 \| \epsilon \|_n \biggl ({\| \epsilon\|_n   \over \lambda} + 2\| \beta^0 \|_1 
\biggr ) $$
$$ = \biggl ( 1 - { 2 \lambda_0 \over \lambda}   \biggr ) \| \epsilon \|_n^2 -
4 \lambda_0 \| \epsilon \|_n \| \beta^0 \|_1 \ge \eta^2 \| \epsilon \|_n^2 . $$
\hfill $\sqcup \mkern -12mu \sqcap$

\begin{theorem} \label{square-rootLasso.theorem} Let $\lambda_0 \| \epsilon \|_n \ge \| \epsilon^T X \|_{\infty} / n $. 
Suppose that for $\eta =\sqrt 2-1 $ one has $\lambda \eta > \lambda_0$ and
$$ \| \beta^0 \|_1 \le {\| \epsilon \|_n  }\biggl ( {   \lambda \eta -   \lambda_0 \over2 \lambda_0 \lambda } \biggr )  . $$
Define  
$$L:= { \lambda \eta + \lambda_0\over \lambda\eta - \lambda_0} . $$
Then
$$ \| X ( \hat \beta - \beta^0 ) \|_n^2  $$  $$ \le \min_{S} \min_{\beta \in \R^p,\  S_{\beta}   =S  } \biggl \{ \| X( \beta - \beta^0) \|_n^2 + \biggl (  {\lambda \eta  + \lambda_0 
 \over  \eta }\biggr )^2  \biggl ( {\lambda \eta  + \lambda_0 \over 2\lambda_0  } \biggr )^2   { \| \epsilon \|_n^2 |S | \over  \hat \phi^2 (L,S) } \biggr \} .  $$

\end{theorem}

{\bf Proof.} The equation $2\eta = (1- \eta^2) $ gives $\eta= \sqrt 2-1$.
Apply Proposition \ref{square-rootLasso.proposition} and Lemma \ref{square-root-bounds.lemma}
and invoke the 
bound
$$\| \hat \epsilon \|_n \le \| \epsilon \|_n + \lambda \| \beta^0 \|_1$$
$$ \le \| \epsilon \|_n \biggl [ 1+ {\lambda \over 4 \lambda_0 } \biggl ( 1- { 2 \lambda_0 \over \lambda }  - \eta^2 \biggr )
\biggr ]=  \| \epsilon \|_n \biggl ( {\lambda \eta + \lambda_0  \over 2 \lambda_0 } \biggr ) . $$
\hfill $\sqcup \mkern -12mu \sqcap$

Using the same arguments we can formulate a bound for the $\ell_1$-estimation error of the
square-root Lasso.

\begin{theorem} \label{square-rootell1.theorem}
Let $\lambda_0 \| \epsilon \|_n \ge \| \epsilon^T X \|_{\infty} /  n$. 
Suppose that for $\eta =\sqrt 2-1 $ one has $\lambda \eta > \lambda_0$ and
$$ \| \beta^0 \|_1 \le {\| \epsilon \|_n  } \biggl ({  \lambda  \eta -  \lambda_0 \over 2 \lambda_0 \lambda }   \biggr ). $$
Let $0 \le \delta < 1$ be arbitrary and define  for $\lambda\eta  >   \lambda_0  $
$$L:=  { \lambda \eta  + \lambda_0 + \delta ( \lambda\eta  - \lambda_0 )\over (1- \delta) (\lambda\eta - \lambda_0) } . $$
Then 
$$2\delta  (\lambda\eta  - \lambda_0) \| \epsilon \|_n \| \hat \beta - \beta^0 \|_1 + \| X ( \hat \beta - \beta^0 ) \|_n^2   $$ 
$$
\le  \min_{S} \biggl \{ \min_{\beta \in \R^p, \ S_{\beta} = S } \biggl [ 2\delta (\lambda\eta  - \lambda_0 ) \| \epsilon \|_n \| \beta - \beta^0 \|_1+
\| X( \beta - \beta^0) \|_n^2   \biggr ]  $$ 
$$
+  \biggl ( { \lambda \eta + \lambda_0 +  \delta (\lambda\eta  -\lambda_0)  
\over  \eta } \biggr )^2 \biggl ( { \lambda \eta  +\lambda_0   \over 2 \lambda_0 } \biggr )^2
{ \| \epsilon \|_n^2 |S | \over  \hat \phi^2 (L,S)} \biggr \} .  $$
\end{theorem}

{\bf Proof.} Combine Proposition \ref{square-rootLasso.proposition} and Lemma \ref{square-root-bounds.lemma},
and invoke the arguments of Theorem \ref{ell1.theorem}.
\hfill $ \sqcup \mkern -12mu \sqcap$

\section{Comparison with scaled Lasso}

Let $\lambda >0$ be a fixed tuning parameter. 
Consider the Lasso with scale parameter $\sigma$
$$\hat \beta (\sigma) := \arg \min_{\beta} \biggl \{ \| Y - X \beta \|_n^2 + 2 \lambda \sigma \| \beta \|_1 \biggr \}, $$
the (scale free) square-root Lasso 
$$ \hat \beta_{\sharp}  := \arg \min_{\beta}  \biggl \{ \| Y - X \beta \|_n +  \lambda  \| \beta \|_1 \biggr \} $$
and the scaled Lasso (\cite{sunzhang11})
$$ ( \hat \beta_{\flat} , \tilde \sigma_{\flat}^2 ) := \arg \min_{\beta , \sigma}  \biggl \{ {  \| Y - X \beta \|_n^2
\over \sigma^2}  +  \log \sigma^2 + {2 \lambda  \| \beta \|_1 \over \sigma}  \biggr \}.
 $$
 Then one easily verifies that
 $$\tilde \sigma_{\flat}^2 = \| Y - X \hat \beta_{\flat} \|_n^2 + \lambda \tilde \sigma_{\flat} \| \hat \beta_{\flat} \|_1 $$
 and that
 $\hat \beta_{\flat} = \hat \beta ( \hat \sigma_{\flat}) $.
 Moreover, if we define
 $$\hat \sigma_{\sharp}^2 :=  \| Y - X \hat \beta_{\sharp} \|_n^2  $$
 we see that $ \hat \beta_{\sharp} = \hat \beta ( \hat \sigma_{\sharp } ) $.
 
Let us write $ {1 \over n}  \times$ the residual sum of squares when using $\sigma$ as
scale parameter as 
 $$\hat \sigma^2  (\sigma)  := \| Y - X \hat \beta ( \sigma ) \|_n^2 .$$
 Moreover, write $ {1 \over n}  \times$ residual sum of squares plus penalty when using  $\sigma$ as scale parameter as
 $$ \tilde \sigma^2 (\sigma) :=  \| Y - X \hat \beta ( \sigma ) \|_n^2 + \lambda \sigma \| \hat
 \beta ( \sigma) \|_1 . $$
 Let furthermore
 $$\tilde \sigma_{\sharp}^2 :=  \| Y - X \hat \beta_{\sharp} \|_n^2 + \lambda \hat \sigma_{\sharp} \| \hat
 \beta_{\sharp} \|_1
 $$
 and
 $$ \hat \sigma_{\flat}^2 := \| Y - X\hat \beta_{\flat} \|_n^2  . $$
 The scaled Lasso includes the penalty in its estimator of $\sigma^2$.
 The square-root Lasso does not explicitly estimate $\sigma^2$.
 In any case, in both versions one may decide to include or not the penalty in an
 estimator of $\sigma^2$. If one does one stays on the conservative side.
 
The square-root Lasso obtains $\hat \sigma_{\sharp}^2 $ as a stable point
of the equation $\hat \sigma_{\sharp}^2 = \hat \sigma^2 ( \hat \sigma_{\sharp})$ and the
scaled Lasso obtains $\tilde \sigma_{\flat}^2$ as a stable point of the equation
 $\tilde \sigma_{\flat}^2=  \tilde \sigma^2 ( \tilde \sigma_{\flat} )$.
 By the mere definition of $\tilde \sigma^2 (\sigma)$ and $\hat \sigma^2 (\sigma)$ we also have
 $\tilde  \sigma_{\sharp}^2 = \tilde \sigma^2 ( \hat \sigma_{\sharp} ) $
 and $\hat \sigma_{\flat}^2 = \hat \sigma^2 ( \tilde \sigma_{\flat} ) $. 
 
 We end this section with a lemma showing the relation between the residual sum of squares with
 penalty and
 the correlation between response and residuals.
 
 \begin{lemma}\label{normalequations.lemma} It holds that
 $$\tilde \sigma^2 (\sigma) = Y^T (Y- X \hat \beta (\sigma)) /n .$$
  \end{lemma}
  
  {\bf Proof.} We have 
  $$ Y^T  (Y- X \hat \beta (\sigma)) /n = \| Y- X \hat \beta (\sigma) \|_n^2 + \hat \beta^T X^T ( Y- X \hat
  \beta (\sigma) ) /n $$
  and by the KKT-conditions (see (\ref{KKT.equation}))
  $$
  \hat \beta^T X^T ( Y- X \hat
  \beta (\sigma) )/n= \lambda \sigma \| \hat \beta (\sigma) \|_1 . $$
  \hfill $\sqcup \mkern -12mu \sqcap $
  
  \chapter{Confidence intervals using the Lasso} 

{\bf Abstract} {\it
We establish asymptotic linearity of a 
de-sparsified Lasso. This implies asymptotic normality under certain conditions and
therefore can be used to construct confidence intervals for parameters of interest.
A similar line of reasoning can be invoked to derive bounds in sup-norm for the Lasso and asymptotic
linearity of de-sparsified estimators of a precision matrix.}

\section{Matrix algebra} 
In this section we show the inverse of a symmetric positive definite matrix $\Sigma_0$ in terms
of projections. 

Let $X_0 \in \R^p$ be a random row-vector with distribution $P$.
We let $\Sigma_0 := E X_0^T X_0$ be the inner-product matrix of $X_0$.
If $E X_0=0$ the matrix $\Sigma_0$ is the covariance matrix of $X_0$.
We assume that $\Sigma_0$ is invertible.
Let $\| \cdot \| $ be the $L_2 (P)$-norm.

For each $j \in \{1 , \ldots , p \}$ we define $X_{-j,0} \gamma_j^0$ as the projection of $X_{j,0}$ on
$X_{-j, 0} := \{ X_{k, 0 } \}_{k \not = j }$. Thus
$$  \gamma_j^0 = \arg \min_{\gamma \in \R^{p-1} } \| X_{j,0} - X_{-j,0} \gamma \| . $$
We further define for all $j$
$$C_{k,j}^0 := \cases {  \ \ \ 1 & $k=j $ \cr -\gamma_{k,j}^0 & $k \not= j $\cr  }  $$
and let $C_0:= ( C_{k,j}^0 ) $. The columns of $C^0$ are written as $C_j^0$, $j=1 , \ldots , p $.
Thus
$$X_{j,0} - X_{-j,0 } \gamma_j^0 = X_0 C_j^0 . $$
We call $X_0 C_j^0$ the {\it anti-projection}
of $X_j$, or the vector of {\it residuals}.
The squared length of the anti-projection or {\it residual
variance} is denoted by $\tau_{j,0}^2:= \| X_0 C_j^0 \|^2 $.  Let $\Theta_j^0 := C_j^0 / \tau_{j,0}^2 $, $j=1 , \ldots , p $ and
$\Theta_0 := ( \Theta_1^0 , \ldots , \Theta_p^0)$.
Then $\Theta_0 = \Sigma_0^{-1} $. 

\section{Notation}
We consider a matrix $X$ with $n$ rows and $p$ columns at write
$\hat \Sigma := X^T X / n $. The columns of $X$ are denoted
by $(X_1 , \ldots , X_p)$. We let for a vector $v \in \R^n$ the normalized Euclidean norm be  $\| v \|_n := \sqrt {v^T v /n}$.
For a real-valued function $f$ we let $\| f \|_n^2 = \sum_{i=1}^n f^2 (X_i)/n$.

We often view the matrix $X$ as being random. We then assume that the rows are i.i.d.\ copies of a
random row vector $X_0$ with distribution $P$ and we write $\Sigma_0 := E X_0^T X_0 =
\EE \hat \Sigma $. We moreover write the $L_2 (P)$-norm as $\| \cdot \|$. For a function $f \in
L_2 (P)$ we have $\| f \|^2 = \EE \| f \|_n^2 $. 

For a matrix $A$ we denote it $\ell_{\infty}$-norm by $\| A \|_{\infty} := \max_{k,j} | A_{k,j} | $. We
define the $\ell_1$-operator norm
$$ \vvert A \vvert_1 := \max_j \sum_k | A_{k,j} |  .$$
For matrices $A$ and $B$ the dual norm inequality is
   $$ \| A B \|_{\infty} \le \|A \|_{\infty} \vvert B \vvert_1 . $$
   
   \vskip .1in
   {\bf Asymptotics} To simplify the exposition we sometimes present asymptotic statements ($n \rightarrow \infty$).
   For a sequence $z_n \in \R$ we write that $z_n= {\mathcal O} (1) $ if $\limsup_{n \rightarrow \infty} |z_n | 
   < \infty$. We write $z_n \asymp 1$ if both $z_n = {\mathcal O} (1)$ and $1/z_n = {\mathcal O} (1)$. 
   We write $z_n = o(1)$ if $\lim_{n \rightarrow \infty} z_n =0 $. 
   
   If $Z$ is random variable which is standard normally distributed we sometimes write $Z=
   {\cal N} (0,1) $.

 \section{A surrogate inverse for $\hat \Sigma$}
Consider a $n \times p$ input matrix $X$ with columns $\{ X_j \}_{j=1}^p $.
We define $X_{-j} := \{ X_k \}_{k \not= j } $, $j=1 , \ldots , p $.
Let $\hat \Sigma = X^T X / n $ be the (normalized) Gram matrix.
We consider for each $j$
 $$\hat \gamma_j (\tau_j ) := \arg \min_{\gamma_j } \biggl \{ \| X_j - X_{-j}  \gamma_j \|_n^2 + 2 \underline{\lambda }\tau_j \| \gamma_j \|_1 \biggr \} $$
 the Lasso for node $j$ on the remaining nodes $X_{-j}$ with tuning parameter
 $\underline
 \lambda$
 and scale parameter $\tau_j$. The reason for introducing a scale parameter here is inspired by the aim to use
 a single tuning parameter $\underline \lambda$ for all $p$ node-wise Lasso's. In the next section
 we will employ the square-root node-wise Lasso which corresponds to a particular choice of
 the scales. As we will see this approach leads to a final scale free result.
 
 Denote the normalized residual sum of squares as $\hat \tau_j^2 ( \tau_j) := \| X_j - X_{-j} \hat \gamma_j \|_n^2 $.
 For the square-root node-wise Lasso the equality $  \hat \tau_j^2 ( \tau_j)= \tau_j^2 \not= 0 $ holds.

 Writing $\tau:= {\rm diag} (\tau_1 , \ldots , \tau_p)$
 we define the matrix $\hat \Theta (\tau)$ as
 $$\hat \Theta_{j,j}  (\tau_j) = 1/ \tilde \tau_j^2  (\tau_j) , \ j \in \{ 1 , \ldots , p \}  , $$
 $$ \hat \Theta_{k,j} (\tau_j ) = - \hat \gamma_{k,j} (\tau_j) / \tilde \tau_j^2 (\tau_j), \ k \not= j \in \{ 1 , \ldots , p \}, $$
 with 
 $$\tilde \tau_j^2 := \| X_j - X_{-j} \hat \gamma_j ( \tau_j \|_n^2 + \underline \lambda \tau_j \| \hat
 \gamma ( \tau_j) \|_1 , \ j \in \{ 1 , \ldots , p \}. $$

Let $e_j$ be the $j$-th unit vector and let $\hat \Theta_j (\tau_j ) $ be the $j$-th column of $\hat \Theta (\tau)$.
The following lemma states that $\hat \Theta (\tau)$ can be viewed as {\it surrogate} inverse for the matrix
$\hat \Sigma$. 
 
 \begin{lemma} \label{surrogate.lemma} We have for all $j$
 $$ \| e_j - \hat \Sigma \hat \Theta_j ( \tau_j) \|_{\infty}  \le \underline \lambda \tau_j / \tilde \tau_j^2 (\tau_j)  $$
 and in fact
 $$  e_{j,j}- \left ( \hat \Sigma \hat \Theta_j ( \tau_j)\right)_j  =0 , $$
 $$\biggl |  e_{k,j} -\left ( \hat \Sigma \hat \Theta_j ( \tau_j)  \right )_k \biggr | 
 \le \underline \lambda \tau_j / \tilde \tau_j^2 (\tau_j) , \ k \not= j . $$
 \end{lemma}
 
 {\bf Proof.} From Lemma \ref{normalequations.lemma}
 $$X_j^T ( X_j - X_{-j} \hat \gamma_j (\tau_j) ) /n= \tilde \tau_j^2  ( \tau_j) $$
 so that
 $$X_j^T X \hat \Theta_j (\tau_j) /n= 1 . $$
 Moreover from the KKT-conditions (see (\ref{KKT.equation})) 
 $$X_{-j}^T ( X_j - X_{-j} \hat \gamma_j (\tau_j) ) / n = \underline \lambda \tau_j \hat z_j
 (\tau_j)  , $$
 where $\| \hat z_j (\tau_j)  \|_{\infty} \le 1 $. We may rewrite this as
 $$ X_{-j}^T X \hat \Theta_j (\tau_j) /n = \underline \lambda \tau_j \hat z_j /\tilde \tau_j^2 (\tau_j)  $$
 giving that $\| X_{-j}^T X \hat \Theta_j (\tau_j) \|_{\infty} /n \le \underline \lambda \tau_j / \tilde \tau_j^2
 (\tau_j ) $.
 \hfill $\sqcup \mkern -12mu \sqcap$

 \begin{lemma} \label{variances.lemma} It holds for all $j$ that
$$  \left ( \hat \Theta (\tau)^T \hat \Sigma \hat \Theta (\tau) \right )_{j,j} = { \hat \tau_j^2 (\tau_j) \over
\tilde \tau_j^4 (\tau_j) } . $$
 \end{lemma}
 
 {\bf Proof.} This is simply rewriting the expressions. We have
 $$ \left ( \hat \Theta (\tau)^T \hat \Sigma \hat \Theta (\tau) \right )_{j,j} =
\hat \Theta_j^T(\tau) \hat \Sigma \hat \Theta_j (\tau) 
 = \| X \hat \Theta_j (\tau_j ) \|_n^2 = \| X \hat C_j (\tau_j) \|_n^2 / \tilde \tau_j^4 (\tau_j) , $$
where $\hat C_j (\tau_j):= \hat \Theta_j (\tau_j) \tilde \tau_j^2 (\tau_j) $.
But
$$ \| X \hat C_j (\tau_j) \|_n^2 = \| X_j - X_{-j} \hat \gamma_j (\tau_j ) \|_n^2 = \hat \tau_j^2 (\tau_j) . $$
\hfill $\sqcup \mkern -12mu \sqcap$

  \section{Asymptotic linearity of the de-sparsified Lasso}
  Let $X$ be an $n \times p$ input matrix and $Y$  an $n$-vector of outputs. 
  We let $f^0:= \EE( Y\vert X)$ be the expectation of $Y$ given $X$ and write the noise as $\epsilon=
  Y- f^0$. We assume $X$ has rank $n$ and let $\beta^0$ be any solution of the equation $f^0 = X \beta^0$

Consider the Lasso with scale parameter $\sigma$
$$\hat \beta (\sigma) := \arg \min_{\beta} \biggl \{ \| Y - X \beta \|_n^2 + \bar \lambda \sigma \| \beta \|_1 \biggr \} .$$
 We define as in \cite{zhang2014confidence} or \cite{van2013asymptotically} the {\it de-sparsified} Lasso
 $$\hat b (\sigma, \tau ) = \hat \beta (\sigma) + \hat \Theta^T (\tau) X^T (Y- X \hat \beta (\sigma)) /n . $$
 For $\hat \Theta^T (\tau)$ we choose the square-root node-wise Lasso $\hat \Theta_{\rm \sharp}$
 which for all $j$ has $\hat  \tau_{j, \sharp}^2$ as stable point of the equation
 $$\hat  \tau_{j, \sharp}^2 = \hat \tau_j^2 ( \hat  \tau_{j, \sharp})  = \| X_j - X_{-j} \hat \gamma_j
 (\hat \tau_{j, \sharp} ) \|_n^2  \not= 0 . $$
 We denote the corresponding de-sparsified Lasso as
 $$\hat b_{\sharp} (\sigma ) = \hat \beta (\sigma) + \hat \Theta_{\sharp}^T X^T (Y- X \hat \beta (\sigma)) /n . $$
 The reason for this choice (and not for instance for the scaled node-wise Lasso or a
 node-wise Lasso with cross-validation) is that the problem becomes scale free.
 There remains however the choice of the tuning parameter $\underline \lambda$. 
 Simulations leads to recommending
 the choice $\underline \lambda= \sqrt {\log p / n }$ (a value which is smaller than the common choice for the
 tuning parameter $\bar \lambda$ for the (standard, square-root or scaled) Lasso for $\beta$). 
 
 We now show that up to a remainder term the estimator $\hat b_{\sharp} (\sigma)$ is
 linear. 
 
\begin{theorem} \label{linear.theorem} For all $j$ and for vectors $v_j \in \R^p$ with $\| v_j \|_n = 1$
(depending on $\hat \Sigma$ and $\underline \lambda$ only\footnote{Hence $v_j$ is a fixed (non-random)
known vector.}) with $\| v_j \|_n = 1$ such that
$${ \tilde \tau_{j , \sharp}^2 \over
\hat \tau_{j, \sharp} } \biggl ( \hat b_{j, \sharp} (\sigma ) - \beta_j^0 \biggr ) = \underbrace{v_j^T \epsilon /n}_{\rm 
linear \ term}  + \underbrace{{\rm rem}_j (\sigma)}_{\rm remainder} $$
where the remainder satisfies
 $$\| {\rm rem} (\sigma) \|_{\infty} \le \underline \lambda \| \hat \beta (\sigma) - \beta^0 \|_1 . $$
\end{theorem}

{\bf Proof.} We have
$$\hat b_{\sharp} (\sigma ) = \hat \Theta_{\sharp}^T X^T\epsilon /n +
\hat \beta (\sigma) - \hat \Theta_{\sharp}^T \hat \Sigma ( \hat \beta (\sigma) -
\beta^0)  $$
so for all $j$
$$\hat b_{j, \sharp} (\sigma ) =  \hat \Theta_{j,\sharp}^T X^T\epsilon /n+
\hat \beta_j (\sigma) -\hat \Theta_{j,\sharp}^T \hat \Sigma ( \hat \beta (\sigma) -
\beta^0)  $$
$$= \beta_j^0 + \hat \Theta_{j,\sharp}^T X^T\epsilon /n+ (e_j^T -  \hat \Theta_{j,\sharp}^T \hat \Sigma) ( \hat \beta (\sigma) -
\beta^0)  . $$
We thus find
$${ \tilde \tau_{j , \sharp}^2 \over
\hat \tau_{j, \sharp} } \biggl ( \hat b_{j, \sharp} (\sigma ) - \beta_j^0 \biggr ) = 
 \underbrace{{ \tilde \tau_{j , \sharp}^2 \over
\hat \tau_{j, \sharp} } 
\hat \Theta_{j,\sharp}^T X^T\epsilon /n }_{:= v_j^T \epsilon / n  } +
\underbrace{ { \tilde \tau_{j , \sharp}^2 \over
\hat \tau_{j, \sharp} } (e_j^T -  \hat \Theta_{j,\sharp}^T \hat \Sigma) ( \hat \beta (\sigma) -
\beta^0)}_{:= {\rm rem}_j (\sigma)} 
   $$
where
$$ v_j := { \tilde \tau_{j , \sharp}^2 \over
\hat \tau_{j, \sharp} } 
X \hat \Theta_{j,\sharp}  $$
and
$$ {\rm rem}_j(\sigma) :=
{ \tilde \tau_{j , \sharp}^2 \over
\hat \tau_{j, \sharp} } (e_j^T -  \hat \Theta_{j,\sharp}^T \hat \Sigma) ( \hat \beta (\sigma) -
\beta^0).$$
Invoking Lemma \ref{variances.lemma} we see that
$$ \hat \Theta_{j, \sharp}^T \hat \Sigma \hat \Theta_{j, \sharp}^T
  = ( \hat \Theta_{\sharp} \hat \Sigma \hat \Theta_{\sharp} )_{j,j}   =
 {\hat \tau_{j,\sharp}^2 \over \tilde \tau_{j, \sharp}^4 }. $$
 Therefore $\| v_j \|_n=1   $.
 Moreover by Lemma \ref{surrogate.lemma}, for all $j$
$$ \|e_j^T -  \hat \Theta_{j, \sharp}^T \hat \Sigma  \|_{\infty} \le \underline \lambda {  \hat \tau_{j, \sharp} \over  
\tilde \tau_{j, \sharp}^2}   $$
and hence $\| {\rm rem}(\sigma) \|_{\infty} \le \underline \lambda \| \hat \beta (\sigma) - \beta^0 \|_1 $.

 \hfill $\sqcup \mkern -12mu \sqcap$

\begin{theorem}  \label{normalconfidence.theorem}
Suppose that $\epsilon \sim {\cal N} (0, \sigma^2 I) $. Then for all $j$
$${  \sqrt n  ( \hat b_{j, \sharp} (\sigma ) - \beta_j^0 )  \over
\sigma  \hat \tau_{j, \sharp}/  \tilde \tau_{j , \sharp}^2 } 
= 
  {\cal N} (0,1)+ \Delta_{j} $$
where $\| \Delta \|_{\infty} \le \sqrt n \underline \lambda 
\| \hat \beta( \sigma) - \beta^0 \|_1 / \sigma $.
\end{theorem} 

{\bf Proof.} This follows immediately from Theorem \ref{linear.theorem}.
\hfill $\sqcup \mkern -12mu \sqcap$

{\bf Asymptotics} Suppose that $\epsilon \sim {\cal N} (0, \sigma^2 I) $, that
$\underline \lambda \asymp \sqrt {\log p / n }$ and that
$\| \hat \beta (\sigma) - \beta^0 \|_1 / \sigma = o_{\PP} ( 1/ \sqrt {\log p } )$. Then
$${  \sqrt n  ( \hat b_{j, \sharp} (\sigma ) - \beta_j^0 )  \over
\sigma \hat \tau_{j, \sharp}/  \tilde \tau_{j , \sharp}^2 } 
=  {\cal N} (0,1) + o_{\PP} (1)   . $$
If the noise is not normally distributed one may explore the possibility of applying a central limit theorem to the linear
term. One should then verify the Lindeberg condition as the asymptotics are for triangular arrays. 

\begin{remark} Assuming that the remainder term $\Delta$ in Theorem \ref{normalconfidence.theorem}
is negligible we can apply its result for the construction of confidence intervals. One then needs
a consistent estimator of $\sigma$. One may for example use a preliminary estimator
$\hat \sigma_{\rm pre}$ for the estimation of $\beta^0$  and then apply the normalized residual sum of squares
$\| Y - X \hat \beta (\hat \sigma_{\rm pre} ) \|_n^2$  as variance estimator for the studentizing step.
Alternatively one may choose the tuning parameter by cross-validation resulting in an
estimator $\hat \beta_{\rm cross}$ and then studentize using $\| Y - X \hat \beta_{\rm cross} \|_n$ as 
estimator of scale.
Another approach would be to apply the square-root Lasso or the scaled Lasso for the estimation of $\beta^0$
and $\sigma$ simultaneously. 
\end{remark}

\begin{remark} The parameter $\beta^0$ is generally not identified as we may take it to be any solution of the equation
  $f^0 = X \beta^0$ (in $\R^n$).  However, we can formulate conditions (see also the next remark)
  depending on the particular
  solution $\beta^0$ 
  such that $\| \hat \beta (\sigma) -
  \beta^0 \|_1 $ converges to zero. Such $\beta^0$ are thus {\it nearly identifiable} and
  Theorem \ref{normalconfidence.theorem} can be used to construct
  confidence intervals for nearly identifiable $\beta^0$ which have $\| \hat \beta (\sigma) -
  \beta^0 \|_1 /\sigma$ converging to zero fast enough.
\end{remark}

\begin{remark}
We note that $\| \hat \beta (\sigma)- \beta^0 \|_1 / \sigma$ can be viewed as a properly scaled $\ell_1$-estimation error.
One may invoke Theorem \ref{ell1.theorem} to bound it. According to this theorem we should choose
$\bar \lambda \sigma  >  \lambda_0 \sigma \ge \| \epsilon^T X \|_{\infty} /n $. We then get 
$${ \| \hat \beta (\sigma) - \beta^0 ) \|_1  \over \sigma}  \le
\min_S \biggl \{ \min_{\beta \in \R^p  , \ S_{\beta } = S } \biggl [ 
{ \| \beta - \beta^0 \|_1  \over  \sigma  } $$
$$ + {1 \over 2 \delta  ( \bar \lambda - \lambda_0) } \biggl ( { \| X (\beta - \beta^0 ) \|_n^2 \over  \sigma^2 }+
\biggl [ \bar \lambda  + \lambda_0 + \delta ( \bar \lambda - \lambda_0 ) \biggr ]^2 
{| S| \over \hat \phi^2 ( L, S) } \biggr ) \biggr \} ,$$
with $L= [\bar \lambda + \lambda_0 + \delta ( \bar \lambda - \lambda_0 ) ]/ [(1- \delta) (\bar
\lambda - \lambda_0 ) ]$. 
In particular we have
$${ \| \hat \beta (\sigma) - \beta^0  \|_1  \over \sigma}  \le{ \biggl [ \bar \lambda  + \lambda_0 + \delta ( \bar \lambda - \lambda_0 ) \biggr ]^2 
  \over 2 \delta  ( \bar \lambda - \lambda_0) } {| S_0 | \over \hat \phi^2 (L, S_0) }  
  .$$
  {\bf Asymptotics} If we take
 take $\bar \lambda \asymp \underline \lambda \asymp \lambda_0 \asymp \sqrt {\log p / n} $ and assume
  $1/\hat \phi (L, S_0) = {\mathcal O} (1)$ then (non-)sparseness $|S_0|$ of small order
  $\sqrt n / \log p $ ensures that
  $\|\Delta \|_{\infty} = o_{\PP} (1)$. In other words the remainder term in the linear approximation
  is negligible  if $\beta^0$ is sufficiently $\ell_0$-sparse. But also more generally, if $\beta^0$ is not
  $\ell_0$-sparse one can still have a small enough remainder term by the above trade-off.
   \end{remark}

  \section{Supremum norm bounds for the Lasso with random design}
  The de-sparsified Lasso deals with the bias of the Lasso. We will now highlight this bias
  term for the case of random design. The  bias for fixed design is similar but the fact
  that we then need to choose a {\it surrogate} inverse somewhat obscures the argument.
  
  Let $(X_0,Y_0)\in \R^{p+1}$ with $X_0$ a $p$-dimensional random row-vector and $Y_0 \in \R$.
 The distribution of $(X_0,Y_0)$ is denoted by $P$ and we let $\| \cdot \|$ be the $L_2(P)$-norm.
 We write $\Sigma_0 := E X_0^T X_0$ and 
  assume that $\Sigma_0$ is invertible.  Let $\Theta_0 := \Sigma_0^{-1}:= ( \Theta_1^0 , \ldots , \Theta_p^0)$.
  Define $\tau_{j,0}^2 := 1/ \Theta_{j,j}^0 $ and $C_j^0 := \Theta_j^0 \tau_{j,0}^2 $,
$j=1 , \ldots , p $.
  Define moreover 
  $$\beta^0 := \arg \min_{\beta \in \R^p } \| Y_0 - X_0 \beta \| . $$
  The noise is denoted by $\epsilon_0 := Y_0 - X_0 \beta^0$ and its variance by $\sigma^2 := \| \epsilon_0 \|^2 $. 

We observe a $n \times (p+1)$ matrix $(X,Y)$ and 
we assume in this section that the rows of $(X,Y)$ are i.i.d. copies of $(X_0,Y_0)$.
Then
$$ Y = X \beta^0 + \epsilon $$
where $\epsilon$ consists of i.i.d.\ copies of $\epsilon_0$. 

We write as usual $\hat \Sigma := X^T X / n$.

We examine the  Lasso
$$ \hat \beta := \arg \min_{\beta \in \R^p} \biggl \{ \| Y - X \beta \|_n^2 +2 \lambda \| \beta \|_1 \biggr \} . $$

\begin{lemma} We have
$$ \| \hat \beta - \beta^0 \|_{\infty} \le \| \Theta_0 X^T \epsilon\|_{\infty} /n +
\underbrace{\vvert \Theta_0 \vvert_1 \biggl ( \| \hat \Sigma - \Sigma_0 \|_{\infty} \| \hat \beta - \beta^0 \|_1 + \lambda \biggr )}_{\rm ``bias"} 
 . $$
 \end{lemma}

{\bf Proof.} By the KKT-conditions (see (\ref{KKT.equation}))
$$ \hat \Sigma ( \hat \beta - \beta^0) + \lambda \hat z = X^T \epsilon/n $$
where
$\hat z_j = {\rm sign} ( \hat \beta_j ) $ if $\hat \beta_j \not= 0 $ and
$\| \hat z \|_{\infty} \le 1 $.  
It follows that
$$C_j^{0T } \hat \Sigma ( \hat \beta - \beta^0) + C_j^{0T} \lambda z = C_j^{0T}X^T \epsilon/n $$
and hence
$$C_j^{0T}  \Sigma_0 ( \hat \beta - \beta^0) + C_j^{0T}  ( \hat \Sigma - \Sigma_0)  ( \hat \beta - \beta^0)
+ C_j^{0T} \lambda \hat z = C_j^{0T}X^T \epsilon/n. $$
But
$$C_j^{0T}  \Sigma_0 ( \hat \beta - \beta^0)= \tau_{j,0}^2 (\hat \beta_j - \beta_j^0) . $$
So we find
$$\hat \beta_j - \beta_j^0 = -\Theta_j^{0T}  ( \hat \Sigma - \Sigma_0)  ( \hat \beta - \beta^0)
-\Theta_j^{0T} \lambda \hat z +  \Theta_j^{0T}X^T \epsilon/n  $$
or
$$ \hat \beta - \beta^0 = - \Theta_0 ( \hat \Sigma - \Sigma_0)  ( \hat \beta - \beta^0) -
\lambda \Theta_0 \hat z + \Theta_0 X^T \epsilon/n . $$
It follows that
$$ \| \hat \beta - \beta^0 \|_{\infty} \le  \| \Theta_0 X^T \epsilon\|_{\infty} /n +
\| \Theta_0 ( \hat \Sigma - \Sigma_0)  ( \hat \beta - \beta^0)\|_{\infty} +
\lambda \| \Theta_0 \hat z \|_{\infty} $$
$$ \le \| \Theta_0 X^T \epsilon\|_{\infty} /n  +  \vvert \Theta_0 \vvert_1 \biggl ( \| \hat \Sigma - \Sigma_0 \|_{\infty} \| \hat \beta - \beta^0 \|_1 + \lambda \biggr )
. $$
\hfill $\sqcup \mkern -12mu \sqcap$ 

{\bf Asympototics} Assume that $\lambda \asymp \sqrt {\log p /n}$ and that
$\| \hat \Sigma - \Sigma_0 \|_{\infty} = {\mathcal O}_{\PP} ( \sqrt {\log p /n})$ and
$\|  X^T \epsilon \|_{\infty}/n = {\mathcal O}_{\PP} ( \sqrt {\log p /n})$. Suppose that
$\| \hat \beta - \beta^0 \|_1$ converges to zero in probability. Then
$$ \| \hat \beta - \beta^0 \|_{\infty} \le   (  \| X^T \epsilon \|_{\infty}/n ) \vvert  \Theta_0 \vvert_1 +
\lambda \vvert \Theta_0 \vvert_1 (1+ o_{\PP} (1) )=\sqrt {\log p /n} \vvert \Theta_0  \vvert_1 (1+ o_{\PP}
(1) ) . $$

 \begin{example} (Equal correlation)
   Let $0\le  \rho < 1 $ and
   
   $$\Sigma_0 :=  \pmatrix{ 1&\rho & \cdots & \rho  \cr
   \rho & 1 & \cdots & \rho \cr 
   \vdots & \vdots & \ddots & \vdots \cr 
   \rho & \rho & \cdots & 1 \cr } = (1- \rho) I + \rho \iota \iota^T, \ \iota := \pmatrix {
    1 \cr 1 \cr \vdots \cr 1 \cr }  .$$
   Then
   $$ \Theta_0 = {1 \over 1- \rho} \left \{ I - { \rho \iota \iota^T \over 1- \rho + p \rho } \right \} $$
   and
   $$\vvert \Theta_0 \vvert_1 ={1 \over 1- \rho}  \biggl \{ { 1 + ( 2p-3) \rho \over 1+ (p-1)\rho}\biggr \}  \le {2 \over 1- \rho} . $$
      \end{example}

   \section{Estimating a precision matrix}
 
 In this section we again let $X$ be an $n \times p $ matrix with rows $\{ X_i \}_{i=1}^n $ 
 being i.i.d.\ copies of a random row-vector $X_0 \in \R^p$.
 Write $\hat \Sigma := X^T X/ n$ and $\Sigma_0 := E X_0^T X_0 = \EE \hat \Sigma$ and
 define $W := \hat \Sigma - \Sigma_0$.
 We assume
 $\Theta_0 := \Sigma_0^{-1} $ exists.

 \subsection{The case $ p $ small}
 Suppose $p$ is small so that $\hat \Sigma$ is invertible for $n$ large enough.
Consider the estimator
$$\hat \Theta := \hat \Sigma^{-1} . $$

\begin{lemma} \label{decomposition.lemma} We have the decomposition
$$ \hat \Theta - \Theta_0  
= - \underbrace{  \Theta_0 W \Theta_0}_{\rm linear \ term  } -  \ {\rm rem}_1 , $$ 
where 
$$ \| {\rm rem}_1 \|_{\infty}  \le \| \Theta_0 W \|_{\infty} \vvert \hat \Theta - \Theta_0 \vvert_1   .$$

\end{lemma} 

{\bf Proof. } We may write
$$ \hat \Theta - \Theta_0 = \hat \Sigma^{-1} - \Sigma_0^{-1} =
\Sigma^{-1} \underbrace{( \Sigma_0 - \hat \Sigma)}_{=-W}  \hat \Sigma^{-1}  $$
$$ = - \Theta_0 W  \hat \Theta  =  - \Theta_0 W   \Theta_0 - \Theta_0 W   (\hat \Theta - \Theta_0) . $$
Thus ${\rm rem}_1=  \Theta_0 W ( \hat \Theta - \Theta_0) $ so that
$ \| {\rm rem}_1 \|_{\infty}  \le \| \Theta_0 W \|_{\infty} \vvert \hat \Theta - \Theta_0 \vvert_1   $. 
\hfill $\sqcup \mkern -12mu \sqcap$

{\bf Asymptotics} Suppose $p$ is fixed and in fact that $X_0$ has a fixed distribution $P$ with finite
fourth moments.
Then $ \| \Theta_0 W \|_{\infty} = {\mathcal O}_{\PP} ( 1/ \sqrt n )$ and 
$\vvert \hat \Theta - \Theta_0 \vvert_1 = o_{\PP} (1)$ and hence
$\| {\rm rem}_1 \|_{\infty} = o_{\PP} (1/ \sqrt n )$. Moreover by the multivariate
central limit theorem $\sqrt n ( \hat \Theta - \Theta_0) $ is asymptotically normally
distributed. 

\subsection{The  square-root node-wise Lasso}
We recall the square-root node-wise Lasso.
For $j=1 , \ldots , p $ we consider the square-root Lasso for the regression of the
$j$-th node on  the other nodes with tuning parameter $\underline \lambda$:
$$\hat \gamma_j := \arg \min_{\gamma_j \in \R^p } 
\biggl \{ \| X_j - X_{-j} \gamma_j \|_n + \underline \lambda  \| \gamma_j \|_1 \biggr \} . $$

Write
$$\hat \tau_j := \| X_j - X_{-j} \hat \gamma_j \|_n , \ \tilde \tau_j^2 = \hat \tau_j^2 +
\underline \lambda \hat \tau_j \| \hat \gamma_j \|_1 $$
and
$$ \hat C_{k,j} := \cases {  \ \ \ 1 & $k=j $ \cr -\hat \gamma_{k,j} & $k \not= j $\cr  } . $$
Then
$$ X_{j} - X_{-j} \hat \gamma_j = X \hat C_j . $$

The KKT-conditions read
 $$ \hat \Sigma  \hat C - \pmatrix { \tilde \tau_1^2 & \cdots & 0 \cr
 \vdots & \ddots & \vdots \cr 0 & \cdots & \tilde \tau_p^2 \cr } = \underline \lambda \hat Z
 \pmatrix { \hat \tau_1 & \cdots & 0 \cr
 \vdots & \ddots & \vdots \cr 0 & \cdots & \hat \tau_p \cr }
 $$
 where $\| \hat Z \|_{\infty} \le 1$ and 
 $$\hat Z_{k,j} = \cases { 0 & $ k=j $ \cr
 {\rm sign} ( \hat \gamma_{k,j} ) & $k \not= j , \ \hat \gamma_{k,j} \not= 0 $ \cr 
 } .$$
 
   Let 
   $$\hat \tau := {\rm diag} ( \hat \tau_1 , \ldots , \hat \tau_p) , \ 
   \tilde \tau := {\rm diag} ( \tilde \tau_1 , \ldots , \tilde \tau_p)$$ 
   and
   $$ \hat \Theta := \hat C \tilde \tau^{-2}  .$$
Then we can rewrite the KKT-conditions as
 $$ \hat \Sigma \hat \Theta  - I = \underline \lambda \hat Z \hat \tau \tilde \tau^{-2} . $$
 
We invert the KKT-conditions for the node-wise Lasso to get the 
{\it de-sparsified node-wise Lasso:}
 $$\hat  \Theta_{\rm de-sparse} 
   :=  \hat \Theta + \hat \Theta^T - \hat \Theta^T \hat \Sigma \hat \Theta . $$

 \begin{lemma} \label{decomposition1.lemma} We have
   $$  \hat \Theta_{\rm de-sparse} - \Theta_0  = - \underbrace{ \Theta_0 W \Theta_0}_{\rm linear \ term}
   - {\rm rem}_1 
   -{\rm rem}_2 $$
   where
   $$\|{\rm rem}_1 \|_{\infty} \le \| \Theta_0 W \|_{\infty} \vvert \hat \Theta - \Theta_0 \vvert_1  , \ 
   \| {\rm rem}_2 \|_{\infty} \le  \underline \lambda \| \tilde \tau^{-1} \|_{\infty} \vvert \hat \Theta - \Theta_0 \vvert_1  .  $$
   \end{lemma}
   
   {\bf Proof.} Write
   $$\hat \Theta_{\rm de-sparse} - \Theta_0 = - \Theta_0 W \Theta_0 -
   \underbrace{\Theta_0 W (\hat \Theta  - \Theta_0)}_{:={\rm rem}_1} - \underbrace{( \hat \Theta - \Theta_0)^T( \hat \Sigma \hat \Theta - I )}_{:= {\rm rem}_2}  $$
   and note that
   $${\rm rem}_2 := ( \hat \Theta - \Theta_0)^T \underbrace{( \hat \Sigma \hat \Theta - I )}_{=
   \lambda \hat Z \hat \tau \tilde \tau^{-2} } $$
   Then use that $\| \hat Z \|_{\infty} \le 1 $ and $\| \hat \tau \tilde \tau^{-1} \|_{\infty} \le 1 $.
   \hfill $\sqcup \mkern -12mu \sqcap$
   
   {\bf Asympotics} We have 
   $\| \Theta_0 W \|_{\infty} = O_{\PP} ( \sqrt {\log p / n } )$ under fourth moment conditions. 
   Let $\phi_0^2$ be the smallest eigenvalue of $\Sigma_0$ and let $\sigma_{\max}^2 :=\max_{j}
   \Sigma_{j,j,0} $. Assume that $1/ \phi_0 = {\mathcal O} (1)$ and $\sigma_{\rm max} = {\mathcal O} (1)$.
   If the data are
   Gaussian\footnote{This can be generalized to sub-Gaussian or bounded data.} we get when choosing $\underline \lambda \asymp \sqrt {\log p / n}$ large enough the result
   $\vvert \hat \Theta - \Theta_0 \vvert_1 = {\mathcal O}_{\PP} ( d_0 \sqrt {\log p / n })$
   where $d_0$ is the maximum degree of the nodes. This can be shown using
   Theorem \ref{square-rootell1.theorem} for all node-wise regressions and checking
   that the result is uniform in $j$. To bound $\| \tilde \tau^{-1} \|_{\infty} $
   we may apply the arguments of Lemma \ref{square-root-bounds.lemma}. The final
   conclusion is asymptotic linearity when
   $d_0 = o (\sqrt n / \log p )$.

  \subsection{The graphical Lasso}
  The graphical Lasso (\cite{friedman2008sparse}) is defined as
 $$\hat \Theta = \arg \min_{\Theta \ {\rm p.s.d.}}\biggl \{ 
 {\rm trace} ( \hat \Sigma \Theta )- \log {\rm det} (\Theta) + 2 \lambda \| \Theta \|_{1, {\rm off}}\biggr \}  , $$
 where $\| \Theta \|_{1, {\rm off}}:= \sum_j \sum_{k \not= j} |\Theta_{j,k} |$ and $\lambda >0$ is a
 tuning parameter. The minimization is carried out over all positive semi-definite (p.s.d.) matrices.

The KKT-conditions are now
$$\hat \Sigma - \hat \Theta^{-1} + \lambda \hat Z =0 $$
where
$$\hat Z_{j,k} = {\rm sign} ( \hat \Theta_{j,k} ) , \ \hat \Theta_{j,k} \not= 0 , \ j \not= k , $$
$$ \hat Z_{j,j}  = 0 , $$
$$ \| \hat Z \|_{\infty} \le 1 . $$

We define the de-sparsified graphical Lasso (\cite{Jankova2014}) as
 $$\hat  \Theta_{\rm de-sparse}  := \hat \Theta + \lambda \hat \Theta \hat Z \hat \Theta $$
    $$ \ \ \ \ \ \ \ \ \ \ \  \ \ \ = 2 \hat \Theta - \hat \Theta \hat \Sigma \hat \Theta . $$

\begin{lemma} \label{decomposition2.lemma} We have
 $$ \hat \Theta_{{\rm de}-{\rm sparse}} - \Theta_0  = - \underbrace{\Theta_0 W \Theta_0}_{\rm linear \ term}  - {\rm rem}_1 -
 {\rm rem}_2 $$
 where
 $$\| {\rm rem}_1 \|_{\infty} \le \| \Theta_0 W \|_{\infty} \vvert \hat \Theta - \Theta_0 \vvert_1 
 , \ \| {\rm rem}_2 \|_{\infty} \le \lambda \vvert  \hat \Theta \vvert_1 \vvert \hat \Theta - \Theta_0 \vvert_1 .$$
 
    \end{lemma}
    
  {\bf Proof.} 
  Write
   $$ \hat \Theta_{{\rm de}-{\rm sparse}} - \Theta_0 = \Theta_0 W \Theta_0 -
   \underbrace{\Theta_0 W (\hat \Theta - \Theta_0)}_{:= {\rm rem}_1}  - 
   \underbrace{(\hat \Theta - \Theta_0) ( \hat \Theta^{-1} -
   \hat \Sigma ) \hat \Theta}_{{\rm rem}_2} .$$
 Using the KKT-conditions we get 
   $$ \| {\rm rem}_2 \|_{\infty} =  \lambda\|  ( \hat \Theta - \Theta_0 ) \hat Z \hat \Theta\|_{\infty} 
 \le \lambda \vvert \hat \Theta- \Theta_0 \vvert_1 
    \| \hat Z \hat \Theta\|_{\infty} $$ 
        $$
    \le  \lambda \vvert \hat \Theta- \Theta_0 \vvert_1  \| \hat Z \|_{\infty} \vvert \hat \Theta \vvert_1
    \le \lambda \vvert \hat \Theta- \Theta_0 \vvert_1  \vvert \hat \Theta \vvert_1 .$$
    \hfill $\sqcup \mkern -12mu \sqcap$

   {\bf Asymptotics} Let $d_0$ be the maximal node degree of $\Theta_0$.
   \cite{ravikumar2011high}  show that for $\lambda \asymp \sqrt {\log p / n } $ large enough
   one has under certain (rather restrictive) conditions
   $$ \vvert \hat \Theta - \Theta_0 \vvert_1 = {\mathcal O}_{\PP} ( \lambda d_0 ) .$$
   This implies $ \vvert \hat \Theta  \vvert_1 = {\mathcal O}_{\PP} ( \sqrt {d_0} )$.
  Hence if in addition $d_0^{3/2} = o( \sqrt n / \log p )$
the de-sparsified estimator is asymptotically linear:
$$ \hat \Theta_{\rm de-sparse} - \Theta_0 = \Theta_0 W \Theta_0 + o_{\PP} (n^{-1/2} ) . $$

  \chapter{Chaining including concentration}

{\bf Abstract} {\it We consider chaining and the more general generic chaining method developed by Talagrand, see e.g.\ \cite{talagrand2005generic}.
This allows one to bound suprema of random processes.
Concentration inequalities are refined probability inequalities, mostly again for
suprema of random processes, see e.g.\ \cite{ledoux2005concentration}, \cite{boucheron2013concentration}.
In this chapter we combine the two. We prove a deviation inequality directly using (generic) chaining.}

\section{Notation}
Let ${\cal T}$ be a countable space and for $i=1 , \ldots , n$ and $t \in {\cal T} $ let $X_i (t)$ be a real-valued
random variable. We write
$X_i = X_i ( \cdot)$, $i=1 , \ldots , n$ and $X(t) := (X_1 (t) , \ldots , X_n (t))^T $,
$ t \in {\cal T} $. When ${\cal T}=\{ 1 , \ldots , p \}$ we will often write $X_i(j) := X_{i,j} $,
$j=1 , \ldots , p $ and let $X$ be the matrix $\{ X_{i,j} \}$.  The columns of $X$ are then $\{ X(j) \}_{j=1}^p $,
where $X(j)= (X_1 (j) , \ldots , X_n (j) )^T$. If confusion between observations $X_i$ and variables
$X(j)$ is not likely we often write $X_j$ instead of $X(j)$. 
We assume that the rows of $X$ - the observations $X_1 , \ldots , X_n$ - are independent.

Let $\epsilon_1 , \ldots , \epsilon_n$ be a Rademacher sequence,
that is $\epsilon_1 , \ldots , \epsilon_n$ are independent and $\PP (\epsilon_i=1)=
\PP ( \epsilon_i = -1) = {1 \over 2}$, $i=1,\ldots , n $. Let $\epsilon$ be the vector
$\epsilon:= ( \epsilon_1 , \ldots , \epsilon_n)^T $. We require $\epsilon$ to be independent
of $\{ X_i (t): \ i=1 , \ldots , n , \ t \in {\cal T} \} $. 

For a vector $v \in \R^n$ we use the notation $\| v \|_{n}^2 := v^T v / n $. 

\section{Hoeffding's inequality}

In this section we examine independent real-valued random variables $X_1 , \ldots , X_n$ (i.e., ${\cal T}$ is a singleton). 
Moreover, $c:= (c_1 , \ldots , c_n)^T$ is a vector of positive constants.

\begin{lemma} \label{extreme.lemma} Fix some $i \in \{ 1 , \ldots , n \}$. Suppose that 
$$ \EE X_i = 0, \ |X_i | \le c_i . $$
Then for any convex function $g$ the mean of $g(X_i)$ can be bounded by the mean
of the extremes in $\pm c_i$: 
$$ \EE g(X_i) \le \EE g ( \epsilon_i c_i) . $$
\end{lemma}

{\bf Proof.} 
For all $0 \le \alpha \le 1$ 
and all $u$ and $v$ we
have
$$g(\alpha u + (1- \alpha ) v )\le \alpha g(u)+(1-\alpha) g(v) . $$
Apply this with $\alpha:= \alpha_i := (c_i - X_i )/(2c_i)$, $u:= u_i := - c_i$ and
$v:= v_i = c_i$, $i=1 , \ldots , n$. Then
$X_i =\alpha_i u_i + (1- \alpha_i) v_i = X_i $, 
so
$$g(X_i)  \le \alpha_i g(-c_i)+ (1- \alpha_i ) g(c_i)$$
and since $\EE \alpha_i =1/2$
$$\EE g(X_i)] \le {1 \over 2} g(-c_i) + {1 \over 2} g(c_i)  = \EE g (\epsilon_i c_i) . $$
\hfill $\sqcup \mkern -12mu \sqcap$

\begin{lemma} \label{exp.lemma} Let $X_1 , \ldots , X_n$ be independent random variables
satisfying $\EE X_i =0 $ and $|X_i | \le c_i$ for all $i=1 , \ldots , n $.
Then for all $\lambda >0$
$$ \EE \exp\biggl  [ \lambda  \sum_{i=1}^n X_i  \biggr] \le \EE \exp \biggl [ \lambda 
\epsilon^T c  \biggr ] 
 . $$
\end{lemma}

{\bf Proof.} Let $\lambda >0$. The map $u \mapsto \exp[\lambda u] $ is convex. The result now follows from 
Lemma \ref{extreme.lemma} and the independence assumptions.

\hfill $\sqcup \mkern -12mu \sqcap$

\begin{lemma}\label{expexp.lemma}
For any $z \in \R$,
$${1 \over 2} \exp [ -z] + {1 \over 2} \exp [z] \le \exp [ z^2 / 2] . $$
\end{lemma}
{\bf Proof.}
For any $z \in \R$
$$ {1 \over 2} \exp [-z]+ {1 \over 2} \exp[z] = \sum_{k \ {\rm even}} {1 \over k!} z^k  = \sum_{k=0}^{\infty} {1 \over (2 k)!} z^{2k} . $$
But for $k \in \{0, 1 , 2 , \ldots \}$ we have
$${1 \over (2 k)!}\le {1 \over k! 2^{k} } $$
so that
$${1 \over 2} \exp [-z]+ {1 \over 2} \exp[z] \le \sum_{k=0}^{\infty} { 1 \over k! } {z^{2k} \over 2^k} =
\exp[z^2 / 2 ] . $$
\hfill $\sqcup \mkern -12mu \sqcap$

\begin{lemma} \label{Rademacherexp.lemma}
$$ \EE \exp \biggl [ \lambda 
\epsilon^T c  \biggr ] 
\le  \exp \biggl [ { n \lambda^2 \| c \|_{n}^2 \over  2 } \biggr  ] . $$
\end{lemma} 

{\bf Proof.} 
By Lemma \ref{expexp.lemma} we have
$$ \EE \exp \biggl [ \lambda 
\epsilon^T c   \biggr ]  = \prod_{i=1}^n \biggl ( {1 \over 2} \exp[-\lambda c_i ]+ 
{1 \over 2} \exp[\lambda c_i ] \biggr ) \le \prod_{i=1}^n \exp\biggl [ {\lambda^2 c_i^2 \over  2} \biggr ] $$ $$=
\exp\biggl [{ n \lambda^2 \| c \|_{n}^2  \over 2} \biggr ] . $$
\hfill $\sqcup \mkern -12mu \sqcap$

\begin{theorem} \label{Hoeffdingmoment.theorem}
Let $X_1 , \ldots , X_n$ be independent random variables
satisfying $\EE X_i =0 $ and $|X_i | \le c_i$ for all $i=1 , \ldots , n $.
Then for all $\lambda >0$
$$ \EE \exp\biggl  [ \lambda  \sum_{i=1}^n X_i  \biggr] \le \exp\biggl [n \lambda^2 \| c \|_n^2 / 2 \biggr ]  $$
and for all $a>0$
$$ \PP \biggl ( {1 \over n} \sum_{i=1}^n X_i \ge \| c \|_n \sqrt {2 a \over n}  \biggr ) \le \exp[-a
 ] . $$

\end{theorem} 

{\bf Proof.} The first result follows from combining Lemmas \ref{exp.lemma} and  \ref{Rademacherexp.lemma}.
For the second result we invoke Chebyshev's inequality: for all $\lambda >0$
$$ \PP \biggl ( {1 \over n} \sum_{i=1}^n X_i \ge \| c \|_n \sqrt {2 a \over n}  \biggr ) \le \exp \biggl [\lambda \sqrt {2a}
- \lambda^2 /2 \biggr ] . $$
Now choose $\lambda = \sqrt {2a} $. \hfill $\sqcup \mkern -12mu \sqcap$

\section{The maximum of $p$ averages}\label{pmaximum.section}

We now consider independent  random row vectors $X_1 , \ldots , X_n$ with values in $\R^p$, that is
${\cal T} := \{ 1 , \ldots, p \}$. 
Let
$$ \biggl \|  {1 \over n} \sum_{i=1}^n X_i \biggr \|_{\infty}:= 
\max_{1 \le j \le p} \biggl | {1 \over n} \sum_{i=1}^n X_i (j) \biggr | .
$$

\begin{lemma} \label{Hoeffdingmax.lemma} Assume that for all $j\in \{ 1 , \ldots , p\} $ 
and all $\lambda >0$
$$ \EE \exp\biggl  [ \lambda  \sum_{i=1}^n X_i (j)   \biggr] \le \exp\biggl [n \lambda^2 / 2 \biggr ] . $$
Then
$$  
\EE \biggl \|{1 \over n } \sum_{i=1}^n X_i \biggr \|_{\infty} \le \sqrt {2 \log (2p) \over n } , $$
and for all $a>0$
$$\PP\biggl ( \biggl \| {1 \over n}  \sum_{i=1}^n X_i\biggr  \|_{\infty}  \ge \
  \sqrt {2 ( \log (2p)+a )  \over n} \biggr ) \le \exp[-a] . $$

\end{lemma}\label{Hoeffdingmax2.lemma}
{\bf Proof.} Let $\lambda >0$ be arbitrary. We have
$$\EE \biggl \|{1 \over n } \sum_{i=1}^n X_i \biggr \|_{\infty} =
{1 \over \lambda} \EE \log \exp\biggl [ {\lambda \over n}  \biggl \| \sum_{i=1}^n X_i \biggr \|_{\infty}  \biggr ] $$
$$ \le {1 \over \lambda } \log \EE  \exp\biggl [ { \lambda \over n } \biggl \| \sum_{i=1}^n X_i \biggr \|_{\infty}  \biggr ] $$ 
$$ \le {1 \over \lambda } \log \biggl (  (2p) \exp[  \lambda^2 / (2n) ] \biggr ) =
{\log (2p) \over \lambda } + {\lambda \over 2n  } . $$
Now choose $\lambda = \sqrt { 2n \log (2p) } $.

For the second result one may use the same argument as in the proof of Theorem 
\ref{Hoeffdingmoment.theorem} to find that for all $a >0$ and all $j$
$$ \PP \biggl ( {1 \over n} \biggl  | \sum_{i=1}^n X_i(j) \biggr | \ge \sqrt { 2a \over n} \biggr ) \le
2 \exp [-a] . $$
Then the inequality for the maximum over $j$ follows from the union bound. 
\hfill $\sqcup \mkern -12mu \sqcap$

\section{Expectations of positive parts} 

We let $[ x ]_+ := x \vee 0 $ denote the positive part of $x \in \R$. 

\begin{lemma} \label{positivepart.lemma} Let for some
$S \in \Nat$ and for $s=1 , \ldots , S $, $Z_s$ be non-negative random variables that satisfy for
certain positive constants $\{ H_s\} $ and for all $a>0$
$$\PP ( Z_s \ge \sqrt {H_s+ a}) \le \exp[-a] . $$
Then
$$\EE \max_{1 \le s \le S} \biggl [ Z_s - \sqrt {H_s +s \log 2 } \biggr ]_+^2 \le 1 . $$

\end{lemma}

{\bf Proof.} Clearly
$$\EE \biggl [ Z_s -  \sqrt {H_s +s \log 2 } \biggr ]_+^2=
\int_0^{\infty} \PP \biggl ( Z_s -  \sqrt {H_s +s \log 2 } \ge \sqrt a \biggr ) da $$
$$\le \int_{0}^{\infty} \PP \biggl ( Z_s \ge \sqrt {H_s + s \log 2  +a} \biggr ) da \le
\int_{0}^{\infty} \exp[-s \log 2 -a ] da = 2^{-s} . $$
Hence
$$\EE \max_{1 \le s \le S} \biggl [ Z_s -  \sqrt {H_s + s \log 2 } \biggr ]_+^2 \le 
\sum_{s=1}^S  \EE\biggl [ Z_s -  \sqrt {H_s + s \log 2 } \biggr ]_+^2 \le \sum_{s=1}^S 2^{-s} \le 1 . $$
\hfill $\sqcup \mkern -12mu \sqcap$ 

\section{Chaining using covering sets}

Consider a subset of a metric space $({\cal T} , d)$.
Let $X_1 , \ldots , X_n \in L_{\infty} ({\cal T}) $ and $X(t) = ( X_1 (t) , \ldots , X_n (t) )^T$, $t \in {\cal T}$. 

Fix some $t_0 \in {\cal T}$ and denote the {\it radius} of ${\cal T}$ by
$$R_n := \sup_{t \in {\cal T}} d( t, t_0)  . $$
Fix some $S \in \Nat$.
For each $s=1 , \ldots , S $ we let ${\cal G}_s \subset {\cal T} $ be a $2^{-s} R_n$ covering set of
${\cal T}$. We take ${\cal G}_0 := \{ t_0\}  $. 
For a given $t$ we let its parent in ${\cal G}_S$ be
$$w(t, S) := \arg \min_{w \in {\cal G}_S } d(t, w) . $$
We let for $s=1, \ldots , S$, $w(t,s-1)$ be the parent of $w(t,s)$:
$$ w(t, s-1):= \arg \min_{w \in {\cal G}_{s-1} } d( w(t,s), w) .$$ Hence $w(t,0) = t_0 $ for all $t$.

\section{Generic chaining}

We let $\{ {\cal G}_s\}_{s \in \Nat} $ be a 
sequence of  finite non-empty subsets of ${\cal T} $ and
we let ${\cal G}_0 = \{ t_0 \} $. 
Consider maps $t \mapsto w(t,s) \in {\cal G}_s $,
 $s=0, \ldots , S$.

We let
$$H_s:= \log (2| {\cal G}_s | ) , \ s =1, \ldots , S , $$ and
$$ \gamma_n  (S) := \max_{j \in {\cal G}_S } \sum_{s=1}^S d_j(s) \sqrt {2(H_s  + s \log 2)/n }$$
where $d_j (s) := d ( w(j,s), w(j,s-1))$ ($j \in {\cal G}_S $, $s=1, \ldots , S $). Define
$$R_n (S):= \max_{j \in {\cal G}_S} \sum_{s=1}^S d_j (s) . $$

Consider the averages $\bar X_n (t) := \sum_{i=1}^n X_i (t)/n$ and 
$\bar \delta_n (t,S) := \bar X_n (t) - \bar X_n( w(t,S))$. 
\begin{theorem}\label{generic-expectation.theorem}
Assume that for each $t , \tilde t \in {\cal T}$ and all $\lambda >0$,
$$\EE \exp \biggl [ \lambda\biggl  | \sum_{i=1}^n ( X_i (t) - X_i (\tilde t)) \biggr  | \biggr ] \le 2 \exp\biggl [ {n \lambda^2 d^2 (t, \tilde t ) \over 2 } \biggr ] . $$
Then
 $$\EE \| \bar X_n -\bar X_n (t_0)  \|_{\infty}  \le  \gamma_n (S)+
 \sqrt {2 R_n^2 (S) / n } + \EE  \| \bar \delta_n ( \cdot , S ) \|_{\infty} . $$
\end{theorem}

{\bf Proof.} 
We may write for all  $t \in {\cal T}$ 
$$ \bar X_n (t) - \bar X_n (t_0) = \sum_{s=1}^{S} \biggl (\bar X_n (w( t,s )) - \bar X_n ( w(t,s-1) )   \biggr ) +  \bar \delta_n (t,S) . $$

We have
$$ \sup_{t \in {\cal T}  }|  \bar X_n (t) - \bar X_n (t_0) | \le
\max_{j \in {\cal G}_S } \sum_{s=1}^{S} \left | \bar X_n (w( j,s )) - \bar X_n ( w(j,s-1) )  \right | + 
\| \bar \delta_n  (\cdot , S\|_{\infty} $$
 But
 $$ \max_{j \in {\cal G}_S } \sum_{s=1}^{S} \left  |  \bar X_n (w(j,s)) -
 \bar X_n ( w(j, s-1) )  \right |   - \gamma_n (S)  $$ $$\le 
 \max_{j {\cal G}_S } \sum_{s=1}^S \biggl ( { |  \bar X_n (w(j,s)) - \bar X_n ( w(j, s-1) ) |  \over
  d_j(s) } - \sqrt {2(H_s+ s \log 2) \over n }  \biggr )  d_j(s)   $$
$$ \le \max_{j \in{\cal G}_S  } \sum_{s=1}^S \biggl [  {|  \bar X_n (w(j,s)) - \bar X_n ( w(j, s-1) ) |  \over
  d_j(s) }
- \sqrt {2(H_s+ s \log 2 )\over n}  \biggr ]_+   d_j(s)  $$
$$ \le  \max_{j \in {\cal G}_S } \max_{s \in \{ 1 , \ldots , S \} } 
\biggl [ { |  \bar X_n (w(j,s)) - \bar X_n ( w(j, s-1) ) | \over
  d_j (s) } - \sqrt {2(H_s+ s \log 2)\over n }  \biggr ]_+   \sum_{s=1}^S d_j (s)  $$
 $$ \le  \max_{k \in{\cal G}_S }\max_{s \in \{ 1 , \ldots , S \} } 
  \biggl [ { | \bar X_n (w(k,s)) - \bar X_n ( w(k, s-1) ) | \over
  d_k(s)  } - \sqrt {2(H_s+ s \log 2) \over n}  \biggr ]_+ R_n (S). $$
  But by Lemma \ref{Hoeffdingmax.lemma}  for all $a>0$ and $s \in \{ 1 , \ldots , S \} $,
  $$\PP \left ( \max_{k \in {\cal G}_S} { |  \bar X_n (w(k,s)) - \bar X_n ( w(k, s-1) ) | \over
  d_k(s) } \ge \sqrt { 2 (H_s + a) \over n } \right )\le \exp[-a] . $$
    Combine this with Lemma \ref{positivepart.lemma} to find that
  $$ \EE \max_{s \in \{ 1 , \ldots , S \} } 
   \max_{k \in{\cal G}_S}
\biggl [ {   |  \bar X_n (w(k,s)) - \bar X_n ( w(k, s-1) ) |   \over
  d_k(s)  } - \sqrt {2(H_s+ s \log 2)\over n  }  \biggr ]_+ \le \sqrt {2 \over n}  $$
  
  It follows that
   $$\EE \left ( \max_{j \in{\cal G}_S} \sum_{s=1}^{S}   \left  |  \bar X_n (w(j,s)) -
 \bar X_n ( w(j, s-1) )  \right |   \right ) 
   \le   \gamma_n (S) + \sqrt {2 R_n^2 (S) / n } $$ and hence
    $$\EE \| \bar X_n - \bar X_n (t_0) \|_{\infty}  \le  \gamma_n (S)+  \sqrt {2 R_n^2 (S) / n } +
    \EE  \| \bar \delta_n  (\cdot , S) \|_{\infty} . $$
   
   \hfill $\sqcup \mkern -12mu \sqcap$
   
   \section{Concentration}
   We use the same notation as in the previous section. 
   
   \begin{theorem} \label{generic-probability.theorem}
   Assume that for each $t , \tilde t \in {\cal T}$ and all $\lambda >0$,
$$\EE \exp \biggl [ \lambda \biggl | \sum_{i=1}^n (X_i (t) - X_i (\tilde t)) \biggr  | \biggr ] \le 2 \exp\biggl [ {n \lambda^2 d^2 (t, \tilde t ) \over 2 }  \biggr ] . $$
Then for all $a >0$
$$ \PP \biggl (\|  \bar X_n - \bar X_n (t_0) \|_{\infty}  \ge  \gamma_n(a,S) +
\| \bar \delta_n ( \cdot , S) \|_{\infty} \biggr ) \le \exp [-a] $$
where
$$\gamma_n(a,S):= \max_{j \in {\cal G}_S } \sum_{s=1}^S d_j (s) \sqrt {2 H_s +
2 (1+s)(1+a)  \over n } . $$
   
   \end{theorem}
   
   {\bf Proof.} 
   Define for $s=1, \ldots , S$
$$\alpha (s)  :=  \sqrt {2 H_s +
2 (1+s)(1+a)  \over n }  . $$
  Using similar arguments as in the proof of Theorem \ref{generic-expectation.theorem} we find
  $$ \sup_{t \in {\cal T}  }|  \bar X_n (t) - \bar X_n (t_0) | \le
\max_{j \in {\cal G}_S } \sum_{s=1}^{S} \left | \bar X_n (w( j,s )) - \bar X_n ( w(j,s-1) )  \right | + 
\| \bar \delta_n  (\cdot , S) \|_{\infty} $$
 and
 $$ \max_{j \in {\cal G}_S } \sum_{s=1}^{S} \left  |  \bar X_n (w(j,s)) -
 \bar X_n ( w(j, s-1) )  \right |   $$
 $$ =  \max_{j \in{\cal G}_S }  \sum_{s=1}^{S}
  { | \bar X_n (w(j,s)) - \bar X_n ( w(j, s-1) ) | \over
  d_j(s) \alpha (s)   }  d_j (s) \alpha(s) $$
  $$ \le \max_{j \in {\cal G}_S }  \max_{1 \le s \le S} { | \bar X_n (w(j,s)) - \bar X_n ( w(j, s-1) ) | \over
  d_j(s) \alpha (s)  } \sum_{s=1}^S d_j (s) \alpha(s) $$
  $$ \le \max_{k \in {\cal G}_S }  \max_{1 \le s \le S} { | \bar X_n (w(k,s)) - \bar X_n ( w(k, s-1) ) | \over
  d_k(s) \alpha (s)  } \max_{j \in {\cal G}_S } \sum_{s=1}^S d_j (s) \alpha (s) . $$
  
Using Lemma \ref{Hoeffdingmax2.lemma}
  $$
  \PP \biggl ( \max_{j \in G_s } {  | \bar X_n ( w(j,s)) - \bar X_n (w(j,s-1))  |  \over
   d_j (s)  }  \ge \alpha (s)  \biggr ) $$ $$ \le  \exp[-(1+s)(1+a) ] , \ 
   s=1 , \ldots , S . 
  $$
  
  Hence
  $$
  \PP \biggl ( \max_{j \in G_s } \max_{1 \le s \le S} {  | \bar X_n ( w(j,s)) - \bar X_n (w(j,s-1))  |  \over
   d_j (s)  }  \ge \alpha(s)  \biggr ) $$ $$\le  \sum_{s=1}^S \exp[-(1+s)(1+a) ] \le \exp[-a] . $$
   \hfill $\sqcup \mkern -12mu \sqcap$
%
%
  
  \begin{remark} We note that
  $$ \gamma_n (a, S) \le \gamma_n (0,S) + R_n (S) \sqrt { 2a/n } . $$  
  The first term does not depend on $a$ and is a bound for the mean
  of $\| \bar X_n - \bar X_n (t_0) \|_{\infty}$. The second term describes the deviation
from this mean. 
  \end{remark}

\bibliographystyle{plainnat}
\bibliography{reference1}
\end{document}